\documentclass[11pt]{article}
\usepackage{times}
\usepackage{mathfont}
\usepackage{graphicx}

\makeatletter
\if@twoside \oddsidemargin 18pt \evensidemargin 56pt \marginparwidth 100pt
\else \oddsidemargin 36pt \evensidemargin 36pt
\fi

\advance\textheight by 6\baselineskip
\def\section{\@startsection {section}{1}{\z@}{-3.5ex plus -1ex minus
 -.2ex}{2.3ex plus .2ex}{\normalsize\bf}}
\def\subsection{\@startsection{subsection}{2}{\z@}{-3.25ex plus -1ex minus
 -.2ex}{1.5ex plus .2ex}{\normalsize\bf}}
\makeatother

\def\max{\mathop{\rm max}}
\def\min{\mathop{\rm min}}

\newcommand{\gthick}{\overline{\theta}}
\newcommand{\thickness}{\theta}
\newcommand{\bookthick}{\mbox{\it bt\/}}
\newcommand{\proof}{\noindent{\bf Proof \hspace{1.6ex}}}
\newcommand{\QED}{\hspace*{\fill} \rule{0.67em}{0.67em}}
\renewcommand{\theequation}{\thesection.\arabic{equation}}
\begin{document}
\newtheorem{theorem}{Theorem}[section]
\newtheorem{conjecture}{Conjecture}[section]
\newtheorem{observation}{Observation}[section]
\newtheorem{proposition}[theorem]{Proposition}
\newtheorem{lemma}[theorem]{Lemma}
\newtheorem{corollary}[theorem]{Corollary}
\newtheorem{fact}[theorem]{Fact}
\bibliographystyle{plain}

\title{Geometric Thickness of Complete Graphs\thanks{A preliminary
version of this paper appeared in the {\it Sixth Symposium on Graph Drawing, GD '98\/}, (Montr\'eal, Canada, August 1998), Springer-Verlag Lecture Notes in Computer Science 1547, 102--110.}}
\author{Michael B. Dillencourt\thanks{%
Supported by NSF Grants CDA-9617349 and CCR-9703572.},
David Eppstein\thanks{%
Supported by NSF Grant CCR-9258355 and matching funds from Xerox Corp.},
Daniel S. Hirschberg\\[0.5em]
Information and Computer Science\\
University of California\\
Irvine, CA 92697-3425, USA.\\[0.5em]
Email: {\tt \{dillenco,eppstein,dan\}@ics.uci.edu}\/} 


\maketitle

\begin{abstract}
We define the {\it geometric thickness\/} of a graph
to be the smallest number of layers such that we can 
draw the graph in the plane with straight-line edges and
assign each edge to a layer so that no two edges on the same layer cross.  
The geometric thickness lies between two previously studied quantities, 
the (graph-theoretical) thickness and the book thickness.
We investigate the geometric thickness of the family of 
complete graphs, $\{K_n\}$.
We show that the geometric thickness of $K_n$ lies between 
$\lceil (n/5.646) + 0.342\rceil$
and $\lceil n/4 \rceil$, and we give
exact values of the geometric thickness of $K_n$ for 
$n\leq 12$ and $n\in\{15,16\}$.
We also consider the geometric thickness of the family of 
complete bipartite graphs.  In particular, we show that, unlike
the case of complete graphs, there are complete bipartite graphs
with arbitrarily large numbers of vertices for which the geometric
thickness coincides with the standard graph-theoretical thickness.
\end{abstract}

\section{Introduction}
\label{sect:intro}

Suppose we wish to display a nonplanar graph on a color terminal in
a way that minimizes the apparent complexity to a user viewing the
graph.  One possible approach would be to use straight-line edges, color
each edge, and require that two intersecting edges have distinct colors.
A natural question then arises: for a given graph, what is the minimum
number of colors required?

Or suppose we wish to print a circuit onto a circuit board, using uninsulated
wires, so that if two wires cross, they must be on different layers, 
and that we wish to minimize the number of layers required.
If we allow each wire to bend arbitrarily, this problem has been
studied previously; indeed, it reduces to the graph-theoretical thickness
of a graph, defined below. 
However, suppose that we wish to further reduce the complexity of the layout 
by restricting the number of bends in each wire.  In particular, if
we do not allow any bends, then the question becomes:
for a given circuit, what is the minimum number of layers required to print
the circuit using straight-line wires?

These two problems motivate the subject of this paper, namely the 
{\it geometric thickness\/} of a graph.  
We define $\gthick(G)$, the geometric thickness
of a graph $G$, to be the smallest value of $k$ such that we can assign
planar point locations to the vertices of $G$, represent each edge of $G$
as a line segment, and assign each edge to one of $k$ layers so that no
two edges on the same layer cross.
This corresponds to the notion of ``real linear thickness'' introduced
by Kainen \cite{Kain73}.
Graphs with geometric thickness~2 (called ``doubly-linear graphs) 
have been studied by Hutchinson {\it et al.\/}~\cite{Hutc95}, where
the connection with certain types of visibility graphs was explored.


A notion related to geometrical thickness is that of 
(graph-theoretical) thickness of a graph, 
$\thickness(G)$, which has been studied extensively 
\cite{Alek78,Bein97,Cimi95,Dean91,Halt91,Jack83a,Mans83a}
and has been defined as the minimum number 
of planar graphs into which a graph can be decomposed.
The key difference between  geometric thickness and graph-theoretical
thickness is that geometric thickness requires that the vertex
placements be consistent at all layers and that straight-line edges
be used, whereas graph-theoretical thickness imposes no 
consistency requirement between layers.

Alternatively, the graph-theoretical thickness can be defined as the minimum
number of planar layers required to embed a graph such that the 
vertex placements agree on all layers but the edges can be arbitrary
curves \cite{Kain73}.  The equivalence of the two definitions follows from the
observation that, given any planar embedding of a graph, the vertex
locations can be reassigned arbitrarily in the plane without altering
the topology of the planar embedding provided we are allowed to bend
the edges at will \cite{Kain73}.  
This observation is easily verified by induction,
moving one vertex at a time.  

The (graph-theoretical) thickness is now known for all complete graphs
\cite{Alek78,Bein67,Bein65,Maye72,Vasak}, and is given by the following
formula:
\begin{equation}
\thickness(K_n) = \left\{ 
\begin{array}{cl}
1, & 1 \leq n \leq 4\\
2, & 5 \leq n \leq 8\\
3, & 9 \leq n \leq 10\\
\left\lceil\frac{n+2}6\right\rceil, & n > 10
\end{array}
\right. 
\label{eq:thicknessval}
\end{equation}

Another notion related to geometric thickness is 
the {\it book thickness\/} of a graph $G$, $\bookthick(G)$, 
defined as follows \cite{Bern79}.
A {\it book with $k$ pages\/}
or a {\it $k$-book\/}, is a line $L$ (called the {\it spine\/})
in 3-space together with $k$ distinct half-planes (called {\it pages\/})
having $L$ as their common boundary.
A {\it $k$-book embedding\/} of $G$ is an embedding of $G$ in a $k$-book
such that each vertex is on the spine, each edge either lies entirely in 
the spine or is a curve lying in a single page, and no two edges intersect 
except at their endpoints.  
The book thickness of $G$ is then the smallest $k$
such that $G$ has a $k$-book embedding.

It is not hard to see that the book thickness of a graph is equivalent
to a restricted version of the geometric thickness where the
vertices are required to form the vertices of a convex $n$-gon.
This is essentially Lemma 2.1, page 321 of \cite{Bern79}.
It follows that
$\thickness(G) \leq \gthick(G) \leq \bookthick(G)$.  It is shown
in \cite{Bern79} that $\bookthick(K_n)=\lceil n/2\rceil$.

In this paper, we focus on the geometric thickness of complete graphs.
In Section~\ref{sect:ub} we provide an upper bound, 
$\gthick(K_n) \leq \lceil n/4\rceil$.
In Section~\ref{sect:lb} we provide a lower bound. 
In particular, we show that 
$\gthick(K_n) \geq \left\lceil \frac{3-\sqrt{7}}2(n+1)\right\rceil
\geq \left\lceil\frac{n+1}{5.646}\right\rceil$. 
This follows from a more precise expression which gives a slightly
better lower bound for certain values of $n$.

These lower and upper bounds do not match in general.  The smallest
values for which they do not match are $n \in\{13,14,15\}$.
For these values of $n$, the upper bound on $\gthick(K_n)$ 
from Section~\ref{sect:ub} is 4, and the lower bound from 
Section~\ref{sect:lb} is 3.
In Section~\ref{sect:k15},  we  resolve one of these three cases
by showing that  $\gthick(K_{15}) = 4$.
For $n=16$ the two bounds match again, but they are distinct for all
larger $n$.

Section~\ref{sect:bip} briefly addresses the geometric thickness of
complete bipartite graphs; we show that 
\[
\left\lceil\frac{ab}{2a+2b-4}\right\rceil
\leq \thickness(K_{a,b})
\leq \gthick(K_{a,b})
\leq \left\lceil\frac{\min(a,b)}2\right\rceil. 
\]
When $a$ is much greater than $b$, the leftmost and rightmost quantities
in the above inequality are equal.  Hence there are complete bipartite
graphs with arbitrarily many vertices for which the standard thickness 
and geometric thickness coincide.
We also show that the bounds on geometric thickness of complete bipartite
graphs given above are not tight, by showing that
$\gthick(K_{6,6})=2$ and  $\gthick(K_{6,8})=3$.  

Section~\ref{sect:concl} contains a table of the lower and upper bounds
on $\gthick(K_n)$ established in this paper for $n\leq 100$ and
lists a few open problems.

\section{Upper Bounds}
\label{sect:ub}
\setcounter{equation}{0}

\begin{theorem}
$\gthick(K_n) \leq \lceil n/4\rceil$.
\label{thm:ub}
\end{theorem}

\proof
Assume that $n$ is a multiple of 4, and let $n=2k$ (so, in particular,
$k$ is even). 
We show that $n$ vertices can be
arranged in two {\it rings\/} of $k$ vertices each, 
an {\it outer ring\/} and an {\it inner ring\/}, so that
$K_n$ can be embedded using only $k/2$ layers and with no edges on the 
same layer crossing.

\begin{figure}[tb]
\begin{center}
\makebox[0pt]{%
\parbox[t]{2.3in}{\begin{center}
\mbox{\includegraphics[height=3in]{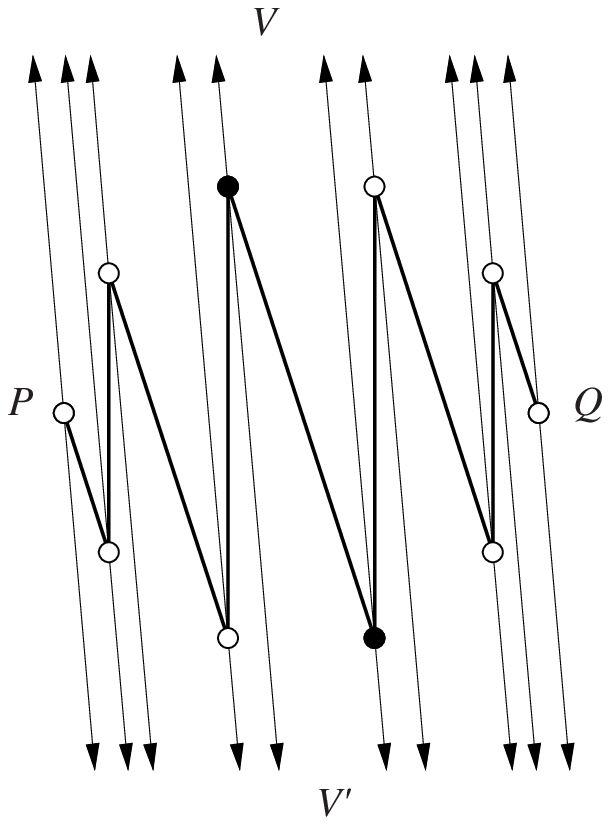}}\\(a)\end{center}}%
\hspace*{0.5in}%
\parbox[t]{2.7in}{\begin{center}
\mbox{\includegraphics[height=3in]{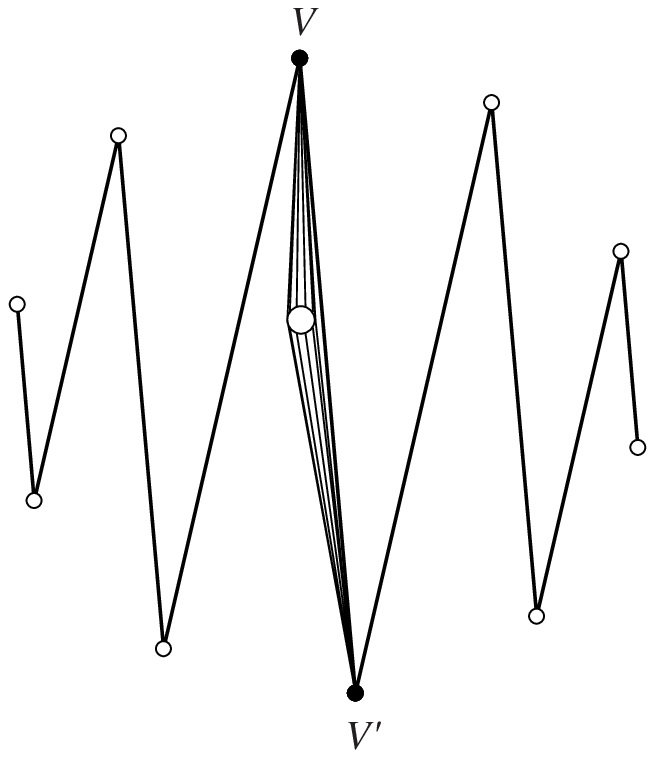}}\\(b)\end{center}}}
\end{center}
\caption[]{Construction for embedding $K_{2k}$ with geometric thickness 
of $k/2$, illustrated for $k=10$. 
(a) The inner ring.
(b) The outer ring.  The circle in the center of (b) represents the inner
ring shown in (a).}
\label{ubfig1}
\noindent\rule{\textwidth}{0.4pt}
\end{figure}

The vertices of the inner ring are arranged to form a regular $k$-gon.
For each pair of diametrically opposite vertices $P$ and $Q$, 
consider the zigzag
path as illustrated by the thicker lines in Figure~\ref{ubfig1}(a).
This path has exactly one diagonal connecting diametrically opposite
points (namely, the diagonal connecting the two dark points in the figure.)
Note that the union of these zigzag paths, taken over all $k/2$ pairs of
diametrically opposite vertices, contains all ${k\choose2}$ edges
connecting vertices on the inner ring.  
Note also that for each choice of diametrically opposite vertices,
parallel rays can be drawn through each vertex, in two opposite directions, 
so that none of the rays crosses any edge of the zigzag path.
These rays are also illustrated in Figure~\ref{ubfig1}(a).

By continuity, if the infinite endpoints of a collection of parallel rays 
(e.g., the family of rays pointing ``upwards'' in 
Figure~\ref{ubfig1}(a)) are replaced by a suitably chosen common
endpoint (so that the rays become segments), the common endpoint can be chosen
so that none of the segments cross any of the edges of the zigzag path.
We do this for each collection of parallel rays, thus forming an outer ring
of $k$ vertices.
This can be done in such a way that the
vertices on the outer ring also form a regular $k$-gon.  By further
stretching the outer ring if necessary, and by moving the
inner ring slightly, the figure can be perturbed so that none of the
diagonals of the polygon comprising the outer ring intersect the
polygon comprising the inner ring.
The outer ring constructed in this fashion is illustrated in
Figure~\ref{ubfig1}(b).

Once the $2k$ vertices have been placed as described above, the
edges of the complete graph can be decomposed into $k/2$ layers.
Each layer consists of:
\begin{enumerate}
\item A zigzag path through the outer ring, as shown
in Figure~\ref{ubfig1}(b).
\item All edges connecting $V$ and $V'$ to vertices of the inner ring, 
where $V$ and $V'$ are the (unique) pair of diametrically opposite points
joined by an edge in the zigzag path through the outer ring.
(These edges are shown as edges connecting the circle with $V$ and $V'$
in Figure~\ref{ubfig1}(b), and as arrows in Figure~\ref{ubfig1}(a)).
\item The zigzag path through the inner ring that does not
intersect any of the edges connecting $V$ and $V'$ with inner-ring
vertices. (These are the heavier lines in Figure~\ref{ubfig1}(a).)
\end{enumerate}

It is straightforward to verify that this is indeed a decomposition
of the edges of $K_n$ into $k/2 = n/4$ layers. \QED

\section{Lower Bounds}
\label{sect:lb}
\setcounter{equation}{0}

\begin{theorem}
\label{thm:thicklb}
For all $n \geq 1$,
\begin{equation}
\gthick(K_n)  \geq \max_{1\leq x \leq n/2}\frac{{n\choose2} - 2{x\choose2}-3}{3n-2x-7}. 
\label{eq:lbminfrac}
\end{equation}
In particular, for $n\geq 12$,
\begin{equation}
\gthick(K_n) \geq \left\lceil \frac{3-\sqrt{7}}2 n +0.342\right\rceil 
\geq \left\lceil\frac{n}{5.646} + 0.342\right\rceil.
\label{eq:lbthickexact}
\end{equation}
\end{theorem}

\proof
We first prove a slightly less precise bound, namely
\begin{equation}
\gthick(K_n) \geq \frac{3-\sqrt{7}}2 n - O(1). \label{eq:lbthickasymp}
\end{equation}
For graph $G$ and vertex set  $X$, let $G[X]$ denote the subgraph of $G$ 
induced by $X$.
Let $S$ be any planar point set, and let $T_1,\ldots T_k$ be a set of 
straight-line planar triangulations of $S$ 
such that every segment connecting two points in $S$
is an edge of at least one of the $T_i$.  
Find two parallel lines that cut $S$ into three subsets $A$, $B$, and $C$
(with $B$ the middle set), with $|A|=|C| = x$, where $x$ is a value to be 
chosen later.  
For any $T_i$, the subgraph $T_i[A]$ is connected,
because any line joining two vertices of $A$ can be retracted onto a path
through $T_i[A]$ by moving it away from the line separating $A$ from $B$.
Similarly, $T_i[C]$ is connected,
and hence each of the subgraphs $T_i[A]$ and $T_i[C]$ has at least $x-1$ edges.  

By Euler's formula, each $T_i$ has at most $3n-6$ edges, so the number of 
edges of $T_i$ not belonging to
$T_i[A]\cup T_i[C]$ is at most 
$3n-6-2(x-1) = 3n-2x-4$.  Hence
\begin{equation}
{n\choose2} \leq 2{x\choose2} + k(3n-2x-4). \label{eq:edgecount1}
\end{equation}
Solving for $k$, we have
\[
k \geq \frac{{n\choose2} - 2{x\choose2}}{3n-2x-4},
\]
and hence
\begin{equation}
k \geq \frac{{n^2} - 2{x^2}}{6n-4x} - O(1). \label{eq:ko1}
\end{equation}
If $x=cn$ for some constant $c$, then the fraction in (\ref{eq:ko1})
is of the form $n(1-2c^2)/(6-4c$).  This is maximized when $c= (3-\sqrt7)/2$.
Substituting the value $x=(3-\sqrt7)n/2$ into (\ref{eq:ko1}) yields
(\ref{eq:lbthickasymp}).

To obtain the sharper conclusion of the theorem, observe that
by choosing the direction of the two parallel lines appropriately, we can
force at least one point of the convex hull of $S$ to lie in $B$.
Hence, of the edges of $T_i$ that  
do not belong to $T_i[A] \cup T_i[C]$,
at least three are on the convex hull.
If we do not count these three edges, then each $T_i$ has at most $3n-2x-7$
edges not belonging to $T_i[A]\cup T_i[C]$,  
and we can strengthen (\ref{eq:edgecount1}) to
\[
{n\choose2} -3\leq 2{x\choose2} + k(3n-2x-7), 
\]
or
\begin{equation}
k \geq \frac{{n\choose2} - 2{x\choose2}-3}{3n-2x-7}. \label{eq:ko2}
\end{equation}
Since 
(\ref{eq:ko2}) holds for any $x$, (\ref{eq:lbminfrac}) follows.

To prove (\ref{eq:lbthickexact}),
let $f(x)$ be the expression on the right-hand side of (\ref{eq:ko2}).
Consider the inequality $f(x) \geq x_0$, where $x_0$ is a constant
to be specified later.
After cross-multiplication, this inequality becomes
\begin{equation}
-x^2+x+\frac{n^2}2 -\frac{n}2 - 3 - (3n - 7 - 2x)x_0 \geq 0. \label{eq:invpar}
\end{equation}
The expression in the left-hand side of (\ref{eq:invpar})
represents an inverted parabola in $x$.
If we let $x=x_0$, we obtain
\begin{equation}
x_0^2 + (8-3n)x_0 + \frac{n^2}2 -\frac{n}2 -3 \geq 0,  \label{eq:defx0}
\end{equation}
and if we let $x=x_0+1$ we obtain the same inequality.
Now, consider $x_0$ of the form $An+B-\epsilon$.
Choose $A$ and $B$ so that if $\epsilon=0$, the terms involving $n^2$
and $n$ vanish in (\ref{eq:defx0}).  This gives the values
$A=(3-\sqrt7)/2$ and $B=\sqrt7(23/14)-4$.
Substituting $x_0=An+B-\epsilon$ with these values of  $A$ and $B$
into~(\ref{eq:defx0}), we obtain
\begin{equation}
\sqrt7\cdot\epsilon\cdot n +(\epsilon^2 -\frac{23\epsilon}{\sqrt7}-3/28)\geq 0. \label{eq:epsdef}
\end{equation}
For $\epsilon=0.0045$, (\ref{eq:epsdef}) will be true when $n \geq 12$.
Therefore, for all $x \in [x_0,x_0+1]$, $f(x) \geq x_0$, when 
$\epsilon=0.0045$ and $n \geq 12$.
In particular, $f(\lceil x_0 \rceil) \geq x_0$.
Since $k$ is an integer,
(\ref{eq:lbthickexact}) follows from (\ref{eq:ko2}).
\QED

\section{The Geometric Thickness of $K_{15}$}
\label{sect:k15}
\setcounter{equation}{0}

The lower bounds on geometric thickness provided by
equation (\ref{eq:lbminfrac}) of Theorem~\ref{thm:thicklb}
are asymptotically larger than the lower bounds on graph-theoretical
thickness provided by equation (\ref{eq:thicknessval}), and they are 
in fact at least as large for all values of $n \geq 12$.
However, they are not tight.
In particular, we show that 
$\gthick(K_{15}) = 4$,
even though (\ref{eq:lbminfrac}) only gives a lower bound of 3.

\begin{theorem}
$\gthick(K_{15}) = 4$.
\label{thm:k15}
\end{theorem}

To prove this theorem, we first note that
the upper bound, $\gthick(K_{15}) \leq 4$, follows immediately from
Theorem~\ref{thm:ub}.

To prove the lower bound, assume that we are given a planar point set
$S$, with $|S|=15$.  We show that there cannot exist a set of three
triangulations of $S$ that cover all ${{15}\choose 2} =105$ line
segments joining pairs of points in $S$.  
We use the following two facts: (1) A planar triangulation with $n$
vertices and $b$ convex hull vertices contains $3n-3-b$ edges; and
(2) Any planar triangulation of a given point set 
necessarily contains all convex hull edges.
There are several cases,
depending on how many points of $S$ lie on the convex hull.

\noindent\underline{Case 1: 3 points on convex hull.}  
Let the convex hull points be $A$, $B$ and $C$.
Let $A_1$ (respectively, $B_1$, $C_1$) be the point furthest from edge 
$BC$  (respectively $AC$, $AB$) within triangle $ABC$.
Let $A_2$ (respectively, $B_2$, $C_2$) be the point next furthest from 
edge $BC$ (respectively $AC$, $AB$) within triangle $ABC$.

\begin{lemma}
The edge $AA_1$ will appear in every triangulation of $S$.
\label{lblemma1}
\end{lemma}

\proof Orient triangle $ABC$ so that edge $BC$ is on the $x$-axis
and point $A$ is above the $x$-axis.
For an edge $PQ$ to intersect $AA_1$, at least one of $P$ or $Q$
must lie above the line parallel to $BC$ that passes through $A_1$.
But there is only one such point, namely  $A$. \QED

\begin{lemma}
At least one of the edges $A_1A_2$ or $AA_2$ will appear in every 
triangulation of $S$.
\label{lblemma2}
\end{lemma}

\proof Orient triangle $ABC$ so that edge $BC$ is on the $x$-axis
and point $A$ is above the $x$-axis.
For an edge $PQ$ to intersect $A_1A_2$ or $AA_2$, at least one of $P$ or $Q$
must lie above the line parallel to $BC$ that passes through $A_2$.
There are only two such points, $A$ and $A_1$.
Hence an edge intersecting $A_1A_2$ must necessarily be $AX$ 
and an edge intersecting $AA_2$ must necessarily be $A_1Y$, 
for some points $X$ and $Y$ outside triangle $AA_1A_2$.
Since edges $AX$ and $A_1Y$ both split triangle $AA_1A_2$,
they intersect, so both edges cannot be present in a triangulation.
It follows that either $A_1A_2$ or $AA_2$ must be present. \QED

Now let $Z$ be the set of 12 edges consisting of the three convex hull
edges and the nine edges $pp_1,pp_2, p_1p_2$ (where $p\in\{A,B,C\}$).
Each triangulation of $S$ contains 39 edges, and since any triangulation
contains all three convex hull edges, it follows from 
Lemmas~\ref{lblemma1} and~\ref{lblemma2} that at least 9 edges of any 
triangulation must belong to $Z$.  Hence a triangulation contains at most
30 edges not in $Z$.  Thus three triangulations can contain at most
$30\cdot3 + 12 = 102$ edges, and hence cannot contain all 105 edges
joining pairs of points in $S$.

\noindent\underline{Case 2: 4 points on convex hull.}  
Let $A$,$B$,$C$,$D$ be the four convex hull vertices.  
Assume triangle $DAB$ has at least one point of $S$ in its interior
(if not, switch $A$ and $C$).
Let $A_1$ be the point inside triangle $DAB$
furthest from the line $DB$.
By Lemma~\ref{lblemma1}, the edge $AA_1$ must appear in every triangulation
of $S$, as must the 4 convex hull edges.  Since any triangulation of $S$
has 38 edges, three triangulations can account for at most $3\cdot33 + 5 =104$
edges.

\noindent\underline{Case 3: 5 or more points on convex hull.}  
Let $h$ be the number of points on the convex hull.  A triangulation
of $S$ will have $42-h$ edges, and all $h$ hull edges must be in each
triangulation.  So the total number of edges in three triangulations
is at most $3(42-2h) + h = 126-5h$, which is at most 101 for $h \geq 5$. 

This completes the proof of Theorem~\ref{thm:k15}.

\section{Geometric Thickness of Complete Bipartite Graphs}
\label{sect:bip}
\setcounter{equation}{0}

In this section we consider the geometric thickness of complete
bipartite graphs, $K_{a,b}$.  
We first give an upper bound, (Theorem~\ref{thm:bip}); it is convenient
to state this bound in conjunction with the obvious lower bound on standard
thickness that follows from Euler's formula.
It follows from this theorem that, for any $b$,
$\gthick(K_{a,b}) = \thickness(K_{a,b})$
provided $a$ is sufficiently large
(Corollary~\ref{cor:bip}).
Hence, unlike the situation with complete graphs,  
there are complete bipartite graphs with arbitrarily many
vertices for which the standard thickness and geometric thickness
coincide.
We show that the lower bound in
Theorem~\ref{thm:bip} is not a tight bound for geometric thickness
by showing that $\gthick(K_{6,8}) = 3$.
A pair of planar drawings demonstrating that $\thickness(K_{6,8}) = 2$
can be found in \cite{Mans83a}.
Finally we show that the upper bound in
Theorem~\ref{thm:bip} is also not tight, since 
$\gthick(K_{6,6}) = 2$ while 
Theorem~\ref{thm:bip} only implies that 
$\gthick(K_{6,6}) \leq 3$.

\begin{theorem}
\label{thm:bip}
For the complete bipartite graph $K_{a,b}$,
\begin{equation}
\left\lceil\frac{ab}{2a+2b-4}\right\rceil
\leq \thickness(K_{a,b})
\leq \gthick(K_{a,b})
\leq \left\lceil\frac{\min(a,b)}2\right\rceil. \label{eq:bip1}
\end{equation}
\end{theorem}

\proof The first inequality follows from Euler's formula, since a planar
bipartite graph with $a+b$ vertices can have at most $2a+2b-4$ edges.
To establish the final inequality, assume that $a \leq b$
and $a$ is even.  
Draw $b$ blue vertices in a horizontal line, with $a/2$ red vertices 
above the line and $a/2$ red vertices below.  
Each layer consists of all edges connecting the
blue vertices with one red vertex from above the line and one
red vertex from below. \QED

\begin{corollary}
\label{cor:bip}
For any integer $b$, 
$\gthick(K_{a,b}) = \thickness(K_{a,b})$ provided
\begin{equation}
a > \left\{ \begin{array}{ll}
             \frac{(b-2)^2}2, &\quad\mbox{if $b$ is even}\\[0.3em]
             (b-1)(b-2), &\quad\mbox{if $b$ is odd}
            \end{array}
   \right.
\end{equation}
\end{corollary}

\proof 
If $a>b$, the leftmost and rightmost quantities 
in (\ref{eq:bip1}) will be equal provided
$ab/(2a+2b-4)>(b-2)/2$ if $b$ is even, or provided
$ab/(2a+2b-4)>(b-1)/2$ if $b$ is odd.
By clearing fractions and simplifying, we see that this happens when
(\theequation) holds.
\QED

\begin{theorem}
\label{thm:k68}
$\gthick(K_{6,8}) = 3$.
\end{theorem}

\proof 
It follows from the second inequality in Theorem~\ref{thm:bip} that
$\gthick(K_{6,8}) \leq 3$, so we need only show
that $\gthick(K_{6,8}) > 2$.
Suppose that we did have an embedding of 
$K_{6,8}$
with geometric thickness~2, with underlying points set $S$.  Since 
$K_{6,8}$ 
has 14 vertices and 48 edges, and since Euler's formula implies
that a planar bipartite graph with 14 vertices has at most 24
edges,  it follows that each layer has exactly 24 edges and
that each face of each layer is a quadrilateral.

Two-color the points of $S$ according to the bipartition 
of $K_{6,8}$.
We claim that there must be at least one red vertex and one blue vertex 
on the convex hull of $S$.  Suppose, to the contrary, that all 
convex hull vertices are the same color (say red).  Then because each
layer is bipartite and because the convex hull contains at least
three vertices, the outer face in either layer would consist of at least
6 vertices (namely the convex hull vertices and three intermediate blue
vertices), which is impossible because each face is bounded by a 
quadrilateral.  The claim implies that one of the layers (say the first) 
must contain a convex hull edge.  
But then this edge could be added to the second layer without 
destroying either planarity or bipartiteness.  Since the second
layer already has 14 vertices and 24 edges, this is impossible.
\QED

\begin{figure}[tb]
\begin{center}
\mbox{\includegraphics{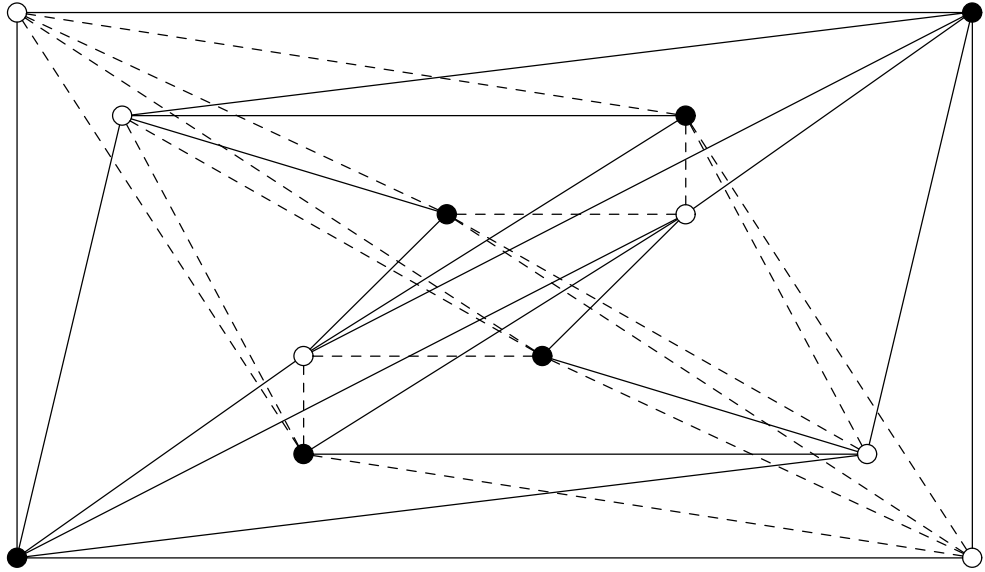}}
\end{center}
\caption[]{A drawing showing that $\gthick(K_{6,6}) = 2$. 
The solid lines represent one layer, the dashed lines the other.}
\label{k66}
\noindent\rule{\textwidth}{0.4pt}
\end{figure}

Figure~\ref{k66} establishes the final claim of the introduction to this 
section, namely that $\gthick(K_{6,6})=2$.

\section{Final Remarks}
\label{sect:concl}
\setcounter{equation}{0}

\begin{table}[bt]

\noindent\rule{\textwidth}{0.4pt}

\caption[]{Upper and lower bounds on $\gthick(K_n)$ established
in this paper.}
\label{tb:numbers}

\begin{center}
\begin{tabular}[t]{|c|r|r|}
\hline
\multicolumn{1}{|c|}{$n$} & \multicolumn{1}{c|}{LB} & \multicolumn{1}{c|}{UB}\\\hline
 1- 4 & 1 &  1\\
 5- 8 & 2 &  2 \\
 9-12 & 3 &  3 \\
13-14 & 3 &  4 \\
15-16 & 4 &  4 \\
17-20 & 4 &  5 \\
21-24 & 5 &  6 \\
25-26 & 5 &  7 \\
27-28 & 6 &  7 \\
29-31 & 6 &  8 \\
  32  & 7 &  8 \\
33-36 & 7 &  9 \\
  37  & 7 & 10 \\\hline
\end{tabular}%
\hspace{0.75in}%
\begin{tabular}[t]{|c|r|r|}
\hline
\multicolumn{1}{|c|}{$n$} & \multicolumn{1}{c|}{LB} & \multicolumn{1}{c|}{UB}\\\hline
38-40 & 8 & 10 \\
41-43 & 8 & 11\\
  44  & 9 & 11\\
45-48 & 9 & 12\\
49-52 &10 & 13\\
53-54 &10 & 14\\
55-56 &11 & 14\\
57-60 &11 & 15\\
61-64 &12 & 16\\
  65  &12 & 17\\
66-68 &13 & 17\\
69-71 &13 & 18\\
  72  &14 & 18\\\hline
\end{tabular}%
\hspace{0.75in}%
\begin{tabular}[t]{|c|r|r|}
\hline
\multicolumn{1}{|c|}{$n$} & \multicolumn{1}{c|}{LB} & \multicolumn{1}{c|}{UB}\\\hline
73-76 &14 & 19\\
  77  &14 & 20\\
78-80 &15 & 20\\
81-82 &15 & 21\\
83-84 &16 & 21\\
85-88 &16 & 22\\
89-92 &17 & 23\\
93-94 &17 & 24\\
95-96 &18 & 24\\
97-99 &18 & 25\\
 100  &19 & 25\\\hline
\end{tabular}
\end{center}

{\small
\noindent
\sloppy\raggedright
{\bf Note:\/}
Upper bounds are from Theorem~\ref{thm:ub}.  The
lower bounds for $n\geq 12$ are from Theorem~\ref{thm:thicklb},
with the exception of the lower bound for $n=15$ which is from
Theorem~\ref{thm:k15}.
Lower bounds for $n < 12$ are from (\ref{eq:thicknessval}).
}

\noindent\rule{\textwidth}{0.4pt}

\end{table}

In this paper we have defined the geometric thickness, $\gthick$, 
of a graph, a measure of approximate planarity that we believe is a natural
notion.  We have established upper bounds and lower bounds on the
geometric thickness of complete graphs.  Table~\ref{tb:numbers} contains
the upper and lower bounds on $\gthick(K_n)$ for $n\leq 100$.

Many open questions remain about geometric thickness.  Here we mention several.
\begin{enumerate}
\item Find exact values for $\gthick(K_n)$ (i.e., remove the
gap between upper and lower bounds in Table~\ref{tb:numbers}).  In particular,
what are the values for $K_{13}$ and $K_{14}$?
\item What is the smallest graph $G$ for which $\gthick(G) > \thickness(G)$?
We note that the existence of a graph $G$ such that $\gthick(G) > \thickness(G)$
(e.g., $K_{15}$) establishes Conjecture 2.4 of \cite{Kain73}.
\item Is it true that $\gthick(G) = O\left(\thickness(G)\right)$ 
for all graphs $G$?
It follows from Theorem~\ref{thm:ub} that this is true for complete graphs.
For the crossing number \cite{Hara62,Pach98}, which like the thickness is a measure of how far
a graph is from being planar, the analogous question is known to have
a negative answer. Bienstock and Dean \cite{Bien93}
have described families of graphs which
have crossing number 4 but arbitrarily high rectilinear crossing number
(where the rectilinear crossing number is the crossing number restricted
to drawings in which all edges are line segments).

\item What is the complexity of computing $\gthick(G)$ for a given graph
$G$?  
Computing $\thickness(G)$ is known
to be NP-complete \cite{Mans83a}, and it certainly seems plausible
to conjecture that the same holds for computing $\gthick(G)$.
Since the proof in \cite{Mans83a} relies heavily on the fact that
$\thickness(K_{6,8}) = 2$,
Theorem~\ref{thm:k68} of this paper shows that this proof cannot
be immediately adapted to geometric thickness.
Bienstock \cite{Bien91} has shown that it is NP-complete to compute 
the rectilinear crossing number of a graph, and that it is NP-hard to 
determine whether the rectilinear crossing number of a given graph
equals the crossing number.
\end{enumerate}

\bibliography{bibliography}
\end{document}